\documentclass[10pt]{amsart}
\usepackage[utf8x]{inputenc}
\usepackage{amsmath, amsthm}
\usepackage{amsfonts}
\usepackage{graphicx,color,amssymb,epstopdf} 
\usepackage{graphicx}
\usepackage{parskip}

\newtheorem{theorem}{Theorem}

\newtheorem{remark}[theorem]{Remark}

\newcommand{\sudda}[1]{}

\newcommand{\F}{{\mathcal F}}
\newcommand{\g}{{\mathfrak g}}

\newcommand{\df}{\mathrm {def}}
\newcommand{\fed}{\mathrm {fed}}

\DeclareMathOperator{\Char}{char}

\title{Software for doing computations in graded lie algebras}
\author{Clas L\"ofwall}\address{Department of Mathematics, Stockholm
  University}\email{clas.lofwall@gmail.com }
  
\author{Samuel Lundqvist}\address{Department of Mathematics, Stockholm
  University}\email{samuel@math.su.se}  \date{}

\begin{document}
\maketitle
\begin{abstract}
We introduce the \emph{Macaulay2} package \emph{GradedLieAlgebras} for doing computations in graded Lie algebras presented by generators and relations.
\end{abstract}

\section{Introduction}
In order to support computer based research on graded Lie (super-)algebras, we have developed the package \emph{GradedLieAlgebras}  as part of \emph{Macaulay2} \cite{M2}.

The package has basic routines for 
 computing Hilbert series, for doing operations on ideals, subalgebras, derivations, and maps. It also has support for constructing holonomy Lie algebras of arrangements, computing homology and constructing minimal models. For a full list of features, we refer to the documentation of the package \cite{package}.

The algorithmic idea used in the package goes back to \cite[Theorem 5.3]{opera}, which was used to identify a periodic structure in a certain 1,2-presented Lie algebra. The first author then 
developed an algorithm and implemented that algorithm in \emph{Mathematica} \cite{mathematica}, under the name \emph{Liedim} \cite{liedim}. That implementation has been cited or referred to in a number of papers, see for instance \cite{fl, llr, peeva, roos},  

The aim of this paper is to describe the  \emph{Macaulay2} implementation, which is a major 
extension of the implementation in \emph{Mathematica}.  

In the next two sections, we discuss implementation details and present the algorithmic theory used in the package. 
In the last section, we give a brief introduction to using the package.

\section{Representing Lie algebras in \emph{Macaulay2}}
In order to be able to use the built-in operations in \emph{Macaulay2}, we decided to convert each computational step in the algorithm to a computation in a corresponding polynomial ring over the same field as the Lie algebra. That polynomial ring is referred to as {\tt lieRing} in the code. 

Let $\g$ be a Lie algebra given by a finite set $\{x_i\}$ of generators and a set of relations. The generators have predefined degrees given by a function $\deg: \{x_i\} \to \mathbb{Z}_+$ and the relations are supposed to be homogeneous with respect to this degree function. This makes $\g$ a \emph{positively graded Lie algebra}. When $\g$ is a Lie super-algebra, the generators have an additional $\mathbb{Z}/2\mathbb{Z}$-grading, and the relations are then supposed to be homogeneous also with respect to this grading.
 
 An iterated Lie product in $\g$ of the form $[x_{i_1},[x_{i_2},[x_{i_3},\ldots [x_{i_{m-1}},x_{i_m}] \ldots ]]]$ is called a \emph{Lie monomial}. 
 In the program, Lie monomials are identified with monomials in {\tt lieRing}, and the Lie product of two elements is performed by a repeated series of normal form computations in {\tt lieRing}, where the normal form computations are being performed with respect to a family of Gr\"obner bases in {\tt lieRing}. It is important to understand that we use the Gr\"obner bases only as a way of doing Gaussian elimination, and that there is no connection to the Buchberger algorithm.

We now describe this correspondence in detail. To make the notation more easy to follow, we will use a slightly different way of naming the generators than in the program.

 Each generator $x_i$ in $\g$ corresponds to $n$ generators $x_{i 1},\ldots,x_{i n}$ in {\tt lieRing}, where $n$ is an upper bound on the degrees handled during our computations. 
 A Lie monomial 
 $[x_{i_1},[x_{i_2},[x_{i_3},\ldots [x_{i_{m-1}},x_{i_m}] \ldots ]]]$ in $\g$ is
 represented as a monomial in {\tt lieRing} in the following way. A generator $x_i$ is represented by $x_{i1}$, and if $e$ is a Lie monomial of degree $d-1$ represented by $m$, then $[x_i,e]$ is represented by $x_{id}\cdot m$, which is also denoted $x_i.m$.

\begin{remark}
 We were informed by J\"orgen Backelin that similar approaches to coding non-commutative monomials as commutative monomials have been considered independently in \cite{gerasimov, scala}.
\end{remark}
 
 From the algorithm described in Section \ref{sec:alg}, it follows that the basis elements of degree $d$ are of the form $[x_i,e_j]$, with $e_j$ a basis element of degree $<d$. 
 
If $e_j$ is a basis element of degree $r-1$ represented by the monomial $m$ in {\tt lieRing}, and $[x_i,e_j]$ is of degree $d$ but not a basis element, then $x_i.m=x_{i r} \cdot m$ will be the leading monomial of a polynomial in a reduced Gr\"obner basis associated to degree $d$. This polynomial then has the form $x_{i r} \cdot m - \sum c_i m_i$, where each $m_i$ corresponds to a Lie monomial in $\g$ that is a basis element in degree $d$, and where each $c_i$ is an element in the underlying field.
 
 The reduced Gr\"obner basis associated to degree $d$ is the degree $d$ part of the reduced Gr\"obner basis of the ideal generated by the set of polynomials of degree $d$ corresponding to the elements that come from the expressions (\ref{1}), (\ref{2}), (\ref{4}) and (\ref{5}) in the next section.
 
\section{Computing a vector space basis of a graded Lie algebra in a given degree} \label{sec:alg}

A Lie (super-)algebra $\g$ over a field $k$ may be specified by giving a positively
	graded finite set
$X$ of generators and a finite set $Y$ of homogenous relations, 
$$
y=\sum_{x\in X}\lambda_{x,y} x+\sum_i[g_{i,y},h_{i,y}],
$$
where $\lambda_{x,y}\in k$ and $g_{i,y},h_{i,y}\in\F(X)$ -- the free Lie algebra on $X$ over $k$. 
Throughout this section we will use the following example to illustrate the steps in the algorithm. 

Let k=$\mathbb{Z}/3\mathbb{Z}$. The set of generators is $X=\{a,b,c\}$ where $a$ is odd, and $b,c$ are even, and $a,b$ have degree $1$, and $c$ has degree $2$. The set of relations is $Y=\{[a,a],[b,[b,a]]-[a,c]\}$. We want to compute a basis in degree $3$, and we will use \hfill
{\tt lieRing}=$\mathbb{Z}/3\mathbb{Z}[a_1,b_1,c_1,a_2,b_2,c_2,a_3,b_3,c_3]$ with $GRevLex$, where $a_i,b_i$ have degree $1$, and $c_i$ has degree $2$ for $i=1,2,3$.

 The program constructs, degree by degree, starting with degree $1$, a graded $\g$-module $M$ that is a subspace of {\tt lieRing}. The $\g$-module 
 operation on $M$ is denoted $g.m$, where $g\in\g$ and $m\in M$. There is also a $\g$-module map $\df:M\to\g$.
Assume $n\ge1$ and everything is done in degree $<n$. This means that $\df: M_{<n}\to \g/\g_{\ge n}$ is an isomorphism 
of $\g$-modules, with an inverse to be denoted by $\fed$.
 If $\deg(g)+\deg(m)\ge n$, 
 $g.m$ is defined to be zero.  
 
 In the example $\{a_1,b_1\}$ is a basis for $M$ in degree $1$ and $\{b_2a_1,c_1\}$ is a basis for $M$ in degree $2$. We also have the reduction rules 
 $a_2a_1\to 0,\ a_2b_1\to-b_2a_1,\ b_2b_1\to0$. Moreover, $\df(a_1)=a,\df(b_1)=b,\df(b_2a_1)=[b,a],\df(c_1)=c$.
 
To construct $M_n$, the first step is to construct a subspace $\hat M_n$ of {\tt lieRing}, with basis 
$$\{x.m\, |\, x \in X, m \in M, \deg(x)+\deg(m)=n \} \cup \{x_{i1} \, | \, x_i\in X, \deg(x_i)=n \}.$$
In a natural way we get an $\F(X)$-module $\hat M=M_{<n}\oplus \hat M_n$ by first defining the action of $X$, and 
then extending this action by the derivation rule. Also $\df$ is defined on $\hat M_n$ by $\df(x.m)=[x,\df(m)]$ and 
$\df(x_{i1})=x_i$.
It follows that $\df$ is surjective in degree $n$. 

In the example, a basis for $\hat M_3$ is
\begin{align}\{a_3b_2a_1,a_3c_1,b_3b_2a_1,b_3c_1,c_2a_1,c_2b_1\}.
\label{0}
\end{align} 
The next step is to divide out by  a subspace $R_n$ of $\hat M_n$ to obtain a $\g$-module. For this reason we compute
\begin{align}R=Y.\hat M.
\label{1}
\end{align} 
It is a subspace of $\hat M_n$, since $M_{<n}$ is a $\g$-module.
 A generating set for $R_n$ is obtained by computing $y.m$ for each relation $y$ and basis element $m\in M$ such that $\deg(y)+\deg(m)=n$. 
 
 In the example, the space $R_3$ is spanned by $[a,a].a_1=2a.(a_2a_1)=0$ and 
$[a,a].b_1=2a.(a_2b_1)=-2a_3b_2a_1$, yielding the reduction rule $a_3b_2a_1\to0$.

We now apply \cite[Theorem 5.3]{opera} to obtain $M_n$ as $\hat M_n$ modulo $R_n$ and the expressions
\begin{align}
&\sum_{x\in X}\lambda_{x,y}m_x+\sum_ig_{i,y}.\fed(h_{i,y})
 \quad\text{for all }y\in Y_n\label{2},\\
&x.m+\epsilon(m,x)\df(m).m_x\label{3},
\end{align}
where $m_x$ is the element in {\tt lieRing} corresponding to $x\in X$, and
where the last expression is computed for all basis elements $m\in M$ and all $\ x\in X$ such that $\deg(x)+\deg(m)=n$. 
Here $\epsilon(m,x)$ is the sign of interchanging the super-elements $m$ and $x$.

In fact, using a linearization idea described in [7], the ”commutative” law (\ref{3}) --- the computationally most heavy part --- does not need to be checked for all elements. Indeed, in characteristic zero,
a basis is computed for the quotient space $\tilde M_n$ of $\hat M_n$ with respect to (\ref{1}), (\ref{2}), and the extra expressions 
\begin{align} 
g_{i,y}.\fed(h_{i,y})+\epsilon(g_{i,y},h_{i,y})h_{i,y}.\fed(g_{i,y})\label{4} \text{ \,\,\, for all } y \in Y_n \text{ and for all $i$.} 
\end{align}
If $\Char(k)>0$ then a basis for the quotient space $\tilde M_n$ of $\hat M_n$ is computed with respect to (\ref{1}), (\ref{2}), (\ref{4}), and also the expressions
\begin{align} 
x.m+\epsilon(m,x)\df(m).m_x\label{5}
\end{align}
coming from (\ref{3}), for which $\deg x$ is a multiple of the characteristic.

Finally, $M_n$ is obtained from $\tilde M_n$ by factoring out (\ref{3}) applied to the 
basis elements $x.m$ of $\tilde M_n$.  

In the example, (\ref{5}) gives nothing, while (\ref{2}) and (\ref{4}) give
\begin{align*}
b.b_2a_1-a.c_1&=b_3b_2a_1 -a_3c_1\ \implies a_3c_1\to b_3b_2a_1,\\
b.b_2a_1+[b,a].b_1&=b_3b_2a_1-b.(b_2a_1)-a.(b_2b_1)=b_3b_2a_1-b_3b_2a_1-0=0, \\
a.c_1+c.a_1&=a_3c_1+c_2a_1\implies c_2a_1\to - a_3c_1\to -b_3b_2a_1.
\end{align*}
Hence, we have the reduction rules $a_3b_2a_1\to0 \text{ (from } (\ref{1})),a_3c_1\to -b_3b_2a_1,c_2a_1\to b_3b_2a_1$ and hence 
$\tilde M_3$ has the 
basis $\{b_3b_2a_1,b_3c_1,c_2b_1\}$. Finally, (\ref{3}) gives the reduction rule 
$c_2b_1\to-b_3c_1$ yielding the basis $\{b_3b_2a_1,b_3c_1\}$ for $M_3$ and $\df(b_3b_2a_1)=[b,[b,a]]$, $\df(b_3c_1)=[b,c]$.

\section{Using the package}
The main introduction to using the package is by means of the tutorials that are part of the documentation \cite{package}. Here we give three small examples of possible computations. 

The most common way to construct a Lie algebra is by means of the constructor { \tt lieAlgebra}. 
In our first example, we construct the free Lie algebra on three even generators, all of degree $1$.

\begin{verbatim}
i2 : L1=lieAlgebra({a,b,c})
o2 : LieAlgebra
i3 : dims(1,6,L1)
o3 = {3, 3, 8, 18, 48, 116}
i4 : basis(2,L1)
o4 =  {(b a), (c a), (c b)}
\end{verbatim}

Here is the example from Section \ref{sec:alg}. 

%%%%%%%%%%

\begin{verbatim}
i5 : L2=lieAlgebra({a,b,c}, Field=>ZZ/3, Signs=>{1,0,0}, 
         Weights=>{1,1,2})/{a a, b b a - a c}
o5 : LieAlgebra
i6 : dims(1,5,L2)
o6 = {2, 2, 2, 3, 5}
i7 :  b c c a	
o7 = (b a b b c) + (b b a b c) + (b b b b b a)
o7 : L2
i8 :  basis(3,L2)
o8 = {(b b a), (b c)}
\end{verbatim}

Let us now give a short example of computing the homology of a Lie algebra.
The generators are odd, $a$ and $b$ have degree $1$, and homological degree $0$, $c$ has degree $2$, and homological degree $1$. The differential is defined by $a,b \mapsto 0, c \mapsto [a,b]$. The  homology can now be obtained using { \tt lieHomology}, {\tt basis} and {\tt dims} (the columns refer to the first degree, and the rows refer to the homological degree). The Lie subalgebras consisting of the cycles and boundaries of the Lie algebra are obtained using {\tt cycles} and {\tt boundaries}.
The underlying field is $\mathbb{Q}$ by default. \break
\begin{verbatim}
i9  : F3=lieAlgebra({a,b,c},Signs=>1,
         Weights=>{{1,0},{1,0},{2,1}},LastWeightHomological=>true)
o9  : LieAlgebra
i10 : L3 = differentialLieAlgebra{0_F3,0_F3,a b}/{a a, b b}
o10 : L3
o10 : LieAlgebra
i11 : H = lieHomology L3
o11 : H
o11 : VectorSpace
i12 : dims(4,H)
o12 = | 2 0 0 0 |
      | 0 0 2 1 |
      | 0 0 0 0 |
      | 0 0 0 0 |
i13 : basis(4,1,H)
o13 : {(b a c)}
i14 : B = boundaries L3
o14 : B
o14 : LieSubAlgebra
i15 : basis(4,1,B)
o15 : {(a b c) + b a c)}
\end{verbatim}

\vspace{5pt}

\section*{Acknowledgement}

The authors want to thank Dan Grayson for his help with many of the implementation issues we faced during the development of the package.

\end{document}